# A differential game with constrained dynamics and viscosity solutions of a related HJB equation[*]

Rami Atar

Paul Dupuis

The Fields Institute
222 College Street
Toronto, Ontario M5T 3J1
Canada

Lefschetz Center for Dynamical Systems
Brown University
Providence, R.I. 02912
USA

March 22, 1999

**Abstract**

This paper considers a formulation of a differential game with constrained dynamics, where one player selects the dynamics and the other selects the applicable cost. When the game is considered on a finite time horizon, its value satisfies an HJB equation with oblique Neumann boundary conditions. The first main result is uniqueness for viscosity solutions to this equation. This uniqueness is applied to obtain the second main result, which is a unique characterization of the value function for a corresponding infinite time problem. The motivation comes from problems associated with queueing networks, where the games appear in several contexts, including a robust approach to network modeling and optimization and risk-sensitive control.

## 1 Introduction and motivation

Queueing and service networks are used as models in many important application areas, including telecommunication, manufacturing, computer networks and vehicle traffic. The control of such networks is difficult and important, and is one of the major challenges facing modern systems theory.

For many problems of this sort it can be difficult to accurately estimate a system model. Service distributions may be complicated and correlated, arrival rates may vary over time, and so on. In such a case one may wish to consider a *robust* approach to control design. The goal of the present paper is to prove foundational results for such an approach, and in particular to develop some analytic machinery for the characterization and construction of value functions for robust formulations of queueing network optimization problems. These results are also useful for analyzing a number of related problems, including risk-sensitive control of queues and design of ordinary (rather than robust) stabilizing controls.

---

[*]This research was supported in part by the National Science Foundation (NSF-DMS-9704426) and the Army Research Office (DAAH04-96-1-0075).

The formulation we use to design robust controls is an analogue of the $H^\infty$ formulation for unconstrained linear and nonlinear systems. See, for example, Appendix B of [3]. We consider a dynamic model whose evolution can be affected by two controls. One control attempts to keep the system in a good operating region (e.g., bounded queue lengths). Such a control may influence the system through service or routing policies. The other control, analogous to the "disturbance" in $H^\infty$ control, determines the value of various system parameters, such as service and arrival rates. This control will attempt to degrade system performance, and the optimization problem is posed as a game.

When formulating a game of this type one must specify the exact form of the dynamics as well as the cost structure. We have two goals in mind regarding the formulation that is used in this paper. The first is to formulate a class of network optimization problems for which explicit solutions may be possible. The second is to choose a cost that will allow one to prove an analogue of the *small gain theorem* from $H^\infty$ control. The small gain theorem, which is a key result in the robust control literature, precisely quantifies the robustness of controls designed using robust criteria against model perturbations. With the first goal in mind, we use the *Skorokhod Problem* to define the system model. This is the simplest dynamical model which preserves the essential features of a controlled network, and in fact the Skorokhod Problem often arises when one considers a law of large numbers approximation to a more detailed (e.g., discrete state) description of a queueing network. A consequence of using this model is that the value function of the differential game satisfies a first order, nonlinear partial differential equation (PDE) with linear Neumann boundary conditions (see Section 2). With both the first and second goals in mind, we have chosen a cost function that is essentially a "weighted" version of the time till the system state reaches the origin. The qualitative properties of unconstrained (e.g., linear system) and constrained (e.g., Skorokhod Problem) models are fundamentally different, and the standard sort of cost structures used in $H^\infty$ control do not yield an analogue of the small gain theorem. The cost we use appears to be the proper replacement in the setting of constrained processes for the familiar quadratic cost used for unconstrained processes. The robust properties that result from the use of this cost are discussed briefly in Section 5, and a full discussion and proof of an analogue of the small gain theorem may be found in [1].

Very little is known regarding robust approaches to the control of queues. To the authors' knowledge, the only other work of this type is due to Ball, Day and Kachroo [2], which is motivated primarily by problems of automobile traffic control. Although the methods of analysis and cost structure of [2] are different from ours, there are many similarities in the system models.

As noted above, the value function for our differential game satisfies (in a viscosity sense) a first order nonlinear PDE with oblique linear derivative boundary conditions. For several reasons, including ease of implementation, we focus on costs that yield controls that are independent of time. Our main interest is therefore in a stationary optimization problem. Just as in the case of unconstrained processes ([17]) this PDE has poor uniqueness properties (even in the viscosity sense). In order to uniquely characterize the value function among all the possible solutions to the PDE, we relate the value function for this problem to the solution of an associated time dependent problem. As shown in Section 4, under appropriate conditions one can construct a solution to the time dependent problem from the solution to the static problem, and conversely. When the viscosity characterization of the unique solution to the time dependent problem is properly interpreted, it provides a analytic characterization of the particular solution to the static PDE that equals the value function for the differential game.

The organization of the paper is as follows. In the next section we discuss the conditions required for uniqueness of the time dependent partial differential equation, and state and prove the comparison principle. Section 3 introduces the Skorokhod Problem and defines a game whose upper value solves the PDE of Section 2. Various special cases of this game have the same value as the games from queueing network optimization problems that are our main motivation, and this relation is made precise in Section 5. In Section 4 we introduce the static game, indicate its connection to the time dependent game, and characterize its value function. We close in Section 5 with examples. The analytic characterizations of Sections 2 and 4 include as special cases other problems besides robust control, and we illustrate this feature by including a problem concerned with existence of stabilizing controls for a network and an optimization problem that arises in risk-sensitive control.

## 2  Formulation of a PDE

As noted in the introduction, the aim of this paper is to establish uniqueness results for viscosity solutions to a class of PDE. These PDE characterize the value functions of certain differential games that appear in problems of robust and risk-sensitive control of queueing networks. It turns out that just as in analogous problems involving unconstrained processes, one cannot prove uniqueness results for these static (time independent) PDE in the standard viscosity sense. However, there are some interesting features of the class of problems we study which allow us to circumvent this difficulty. In particular, we can define a finite time problem for which there is uniqueness, and under an additional condition obtain a simple set of necessary and sufficient conditions that must be satisfied by the value function for the time independent game.

We first define the class of time dependent PDE. To simplify the exposition we use assumptions that are much stronger than necessary and consider as candidate solutions only Lipschitz continuous functions. The intended use of these results is in the construction of solutions, and for the problems we have in mind the value functions will be Lipschitz continuous.

The domain of interest is the set $\mathbb{R}^d_+ \doteq \{x \in \mathbb{R}^d : x_i \geq 0, i = 1, ..., d\}$, where $d \in \mathbb{N}$ is fixed. For each $i \in \{1, ..., d\}$ we associate with the face $F_i \doteq \{x \in \mathbb{R}^d_+ : x_i = 0\}$ of $\mathbb{R}^d_+$ a vector $\gamma_i$. It is assumed that $(\gamma_i)_i > 0$, and we further assume (without loss) that $(\gamma_i)_i = 1$. There is also an index set $\{1, ..., J\}$, which in the control problem will correspond to the different possible "pure" service assignments (cf. the examples in Section 5).

We next state the time dependent PDE and related boundary and terminal conditions. For notational simplicity we consider only the time interval $[0,1]$. For a function $V : [0,1] \times \mathbb{R}^d \to \mathbb{R}$ with independent variables $(t,x) \in [0,1] \times \mathbb{R}^d$, $DV(t,x)$ denotes the gradient with respect to the spatial variables $x$, and $V_t(t,x)$ denotes the partial derivative with respect to $t$. For reals $a$ and $b$ $a \wedge b$ denotes the smaller of $a$ and $b$ and $a \vee b$ denotes the larger. Let

$$H(p) \doteq \wedge_{j=1}^{J} H_j(p). \tag{2.1}$$

The function $H$ defined in the last display will be the Hamiltonian for the PDE on the interior of $\mathbb{R}^d_+$. Boundary conditions will also be needed, and these turn out to be linear oblique derivative conditions defined in terms of the vectors $\gamma_i$. Let $e_i$ denote the unit vector in coordinate direction $i$ for $i = 1, ..., d$.

We use the following regularity conditions.

**Condition 2.1** *For each $j \in \{1, ..., J\}$ the function $H_j : \mathbb{R}^d \to \mathbb{R}$ is convex. The functions $f : \mathbb{R}^d_+ \to \mathbb{R}$ and $g : [0,1] \times \mathbb{R}^d_+ \to \mathbb{R}$ are bounded and continuous.*

Note that Condition 2.1 implies $H_j(p)$ is finite for all $j \in \{1, ..., J\}$ and $p \in \mathbb{R}^d$.

For a point $x \in \mathbb{R}^d_+$ let $I(x) \doteq \{i : x_i = 0\}$. We can now define (in formal terms) the PDE of interest:

$$\max\{-V_t(t,x) - H(DV(t,x)); V(t,x) - g(t,x)\} = 0, \qquad t \in (0,1), x \in (\mathbb{R}^d_+)^\circ, \tag{2.2}$$

$$\langle DV(t,x), \gamma_i \rangle = 0, \qquad t \in (0,1), i \in I(x), x \in \partial \mathbb{R}^d_+ \tag{2.3}$$

$$V(t,0) = 0 \text{ if } t \in [0,1], \tag{2.4}$$

$$V(1,x) = f(x) \text{ for } x \in \mathbb{R}^d_+ \setminus \{0\}. \tag{2.5}$$

The maximization in equation (2.2) is due to the fact that in the corresponding game we allow the minimizing player to stop the game and pay a stopping cost of $g$. This stopping cost can depend on the time and state of the controlled process.

Let a bounded Lipschitz continuous function $w : [0,1] \times \mathbb{R}^d_+ \to \mathbb{R}$ be given. We say that $w$ is a *viscosity subsolution* [resp., *viscosity supersolution*] to (2.2)-(2.5) if the following conditions hold. Let $\theta : (0,1) \times \mathbb{R}^d_+ \to \mathbb{R}$ be a continuously differentiable function, and let $(s,y) \in (0,1) \times \left(\mathbb{R}^d_+ \setminus \{0\}\right)$ be a local maximum [resp., local minimum] of $w - \theta$. Then

$$\max\left\{\left([-\theta_t(s,y) - H(D\theta(s,y))] \wedge \min_{i \in I(y)} [-\langle D\theta(s,y), \gamma_i \rangle]\right); w(s,y) - g(s,y)\right\} \leq 0 \tag{2.6}$$

[resp.,

$$\max\left\{\left([-\theta_t(s,y) - H(D\theta(s,y))] \vee \max_{i \in I(y)} [-\langle D\theta(s,y), \gamma_i \rangle]\right); w(s,y) - g(s,y)\right\} \geq 0], \tag{2.7}$$

$$w(t,0) = 0 \qquad \text{if } t \in [0,1], \tag{2.8}$$

$$w(1,x) \leq f(x) \qquad [\text{resp.}, w(1,x) \geq f(x)] \qquad \text{if } x \in \mathbb{R}^d_+ \setminus \{0\}. \tag{2.9}$$

If $w$ is both a viscosity subsolution and supersolution, then it is called a viscosity solution.

Now of course we will need some conditions on the vectors $\gamma_i$ that appear in the oblique boundary conditions. We will use the following two conditions. The first condition is intimately connected with the uniqueness theory for the underlying controlled processes (cf. [7, 11] and Section 3). We refer to [11] for detailed information on when Condition 2.2 is satisfied. A set $B \subset \mathbb{R}^d$ is called *symmetric* if $x \in B \Rightarrow -x \in B$.

**Condition 2.2** *Let vectors $\gamma_i, i \in \{1,...,d\}$ that satisfy $\langle \gamma_i, e_i \rangle = 1$ be given. There exists a compact, convex, and symmetric set $B$ with $0 \in B^\circ$, such that if $z \in \partial B$ and if $n$ is an outward normal to $B$ at $z$, then for all $i \in \{1,...,d\}$*

$$\langle z, e_i \rangle \geq -1 \text{ implies } \langle \gamma_i, n \rangle \geq 0 \text{ and } \langle z, e_i \rangle \leq 1 \text{ implies } \langle \gamma_i, n \rangle \leq 0.$$

**Condition 2.3** *The vectors $\gamma_i, i \in \{1,...,d\}$ are linearly independent.*

**Remark** Condition 2.2 is equivalent to the main condition of [7, 11], where it is used to study the Skorokhod Problem that is related to these boundary conditions (cf. Section 3). As phrased in these references, the condition requires that if $z \in \partial B$ and $|\langle z, e_i \rangle| \leq 1$ then $\langle \gamma_i, n \rangle = 0$ for all outward normals $n$ to $B$ at $z$. Clearly Condition 2.2 implies this latter condition. To prove the reverse implication, it suffices to show that $\langle z, e_i \rangle < 0$ implies $\langle \gamma_i, n \rangle \leq 0$. We use the fact that a set which satisfies the condition of [11] is automatically invariant under the oblique projection

$$Lx \doteq x - \langle x, e_i \rangle \gamma_i$$

[11, Section 2.5 of part I]. Thus $y \doteq Lz = z - \langle z, e_i \rangle \gamma_i \in B$. Since $n$ is an outward normal to $B$ at $z$, this implies $-\langle z, e_i \rangle \langle n, \gamma_i \rangle = \langle y - z, n \rangle \leq 0$. Hence if $\langle z, e_i \rangle < 0$ then $\langle \gamma_i, n \rangle \leq 0$, as required.

As is well known, in order to prove a comparison principle for viscosity solutions we need a suitable family of test functions. Our construction of the test functions is similar to that of [8], though much simpler. The test functions are constructed from the set $B$. As a preliminary step, we note that if $B$ satisfies Condition 2.2 then so does the convex set

$$\hat{B} = \frac{1}{1-\delta} \{x : \min\{\|x - y\| : y \in B\} \leq \delta\}$$

for small $\delta > 0$. However, this set always has a smooth boundary, in the sense that if $\hat{n}(z)$ is the set of outward unit normals to $\hat{B}$ at $z \in \partial \hat{B}$, then $\hat{n}(z)$ is single valued, and if $z_i \in \partial \hat{B} \to z \in \partial \hat{B}$, then with an abuse of notation $\hat{n}(z_i) \to \hat{n}(z)$. We can therefore assume without loss of generality that $B$ has a smooth boundary, and we do so for the remainder of the paper.

The test function is constructed so that its level sets are all multiples of $B$. This will force the gradient of the function to inherit certain properties of the outward normals of $B$. In addition, we will need the function to grow in a quadratic fashion.

**Lemma 2.1** *Assume Condition 2.2, and for $c \in [0, \infty)$ let*

$$\Xi(x) = c \Leftrightarrow x \in \partial(cB).$$

*Then there exist $0 < m_1 < M_1 < \infty$ such that if $x \neq 0$ and $y = \alpha x \in \partial B$ for some $\alpha \in (0, \infty)$, then $D\Xi(x) = bn(y)$ for some $b \in [m_1, M_1]$ (where $b$ can depend on $x$ and $n(y)$ is the outward normal to $B$ at $y$). Define the continuously differentiable function $\xi(x) \doteq (\Xi(x))^2$. Then there exist $0 < m < M < \infty$ and a function $b : \mathbb{R}^d \to [m, M]$ such that $D\xi(x) = b(x)\Xi(x)n(y)$, and such that $m\|x\|^2 \leq \xi(x) \leq M\|x\|^2$.*

**Proof** Since $0 \in B^\circ$ and $B$ is bounded and convex, $\Xi$ is well defined. The properties of $\xi$ and $\Xi$ then follow directly from the definitions. ∎

Note that $\Xi$ has linear growth along each radial direction. If we extend the definition of $n$ by letting $\bar{n}(x) = n(y)$ if $x = \alpha y$, $\alpha > 0$, and $y \in \partial B$, then the expression $D\xi(x) = b(x)\Xi(x)\bar{n}(x)$ is valid for all $x \in \mathbb{R}^d$.

In addition to the test function, we also need the following.

**Lemma 2.2** *Assume Conditions 2.2 and 2.3. Then there exists a continuously differentiable function $\mu : \mathbb{R}^d_+ \to [0, 1]$ such that*

$$\langle D\mu(x), \gamma_i \rangle < 0$$

*for all $x \in \mathbb{R}^d_+$ and all $i \in I(x)$, and such that $\sup_{x \in \mathbb{R}^d_+} \|D\mu(x)\| < \infty$.*

**Proof** This lemma would be a special case of [8, Lemma 3.2], except that the lemma in [8] assumes a bounded domain. (Conditions (B.1), (B.4) and (B.5) of [8, Lemma 3.2] hold trivially in our setup, and Conditions (B.6) and (B.8) are equivalent in our setup to Conditions 2.3 and 2.2, respectively.) However, given any compact set $K \subset \mathbb{R}^d_+$, the proof of [8, Lemma 3.2] directly implies the existence of a continuously differentiable function $\nu : \mathbb{R}^d_+ \to [0, \infty)$ with the following properties:

- $\nu$ has compact support,
- Given any $x \in K$, $\langle D\nu(x), \gamma_i \rangle < 0$ for all $i \in I(x)$,
- Given any $x \in \mathbb{R}^d_+$, $\langle D\nu(x), \gamma_i \rangle \leq 0$ for all $i \in I(x)$.

We can guarantee that the range of $\nu$ is contained in $[0, 1]$ by multiplying $\nu$ by an appropriate constant. Hence if $\delta_i \in \mathbb{R}, i \in \mathbb{N}$ satisfy $\delta_i > 0$ and $\sum_{i=1}^\infty \delta_i \leq 1$, then

$$\mu(x) \doteq \sum_{i=1}^\infty \delta_i \nu(\delta_i x)$$

satisfies all the conclusions of Lemma 2.2. ∎

**Theorem 2.3** *Assume Conditions 2.1, 2.2, 2.3, and let $u$ (resp., $v$) be a viscosity subsolution (resp., supersolution) to (2.2)–(2.5). Then $u \leq v$ on $[0, 1] \times \mathbb{R}^d_+$.*

**Remark** Since a viscosity solution is by definition both a sub and a supersolution, the comparison principle stated in Theorem 2.3 implies uniqueness of viscosity solutions.

**Proof** The general layout of the proof follows the standard arguments [6], while the treatment of the boundary conditions uses ideas from the proofs of [8] and [9]. Fix $\tau \in (0,1)$, and for $\delta > 0$ and $c > 0$ let

$$U(t,x) = u(t,x) - \delta \frac{(1-\tau)^2}{(t-\tau)} - c\mu(x),$$

$$V(t,x) = v(t,x) + \delta \frac{(1-\tau)^2}{(t-\tau)} + c\mu(x),$$

where $\mu$ is the function described in Lemma 2.2. To prove the theorem it suffices to show that for all $\delta \in (0,1)$ there is $c_0 \in (0,\infty)$ such that

$$U(t,x) \leq V(t,x) \tag{2.10}$$

if $c \in (0, c_0)$ and $(t,x) \in (\tau, 1] \times \mathbb{R}_+^d$. Let $L < \infty$ denote the larger of the Lipschitz constants of $u$ and $v$. Given $\delta > 0$ we can find $\omega \in (0,1)$ such that

$$|H(p+v_1) - H(p+v_2)| \leq \delta/2 \text{ if } \|p\| \leq L+1 \text{ and } \|v_1\| \vee \|v_2\| \leq \omega. \tag{2.11}$$

Let $\bar{M} = \sup_{x \in \mathbb{R}_+^d} \|D\mu(x)\|$. We then let $c_0 \doteq (\delta/4) \wedge (\omega/4\bar{M})$, and assume henceforth that $c \in (0, c_0)$.

We will use the usual argument by contradiction, and so assume that

$$\rho \doteq \sup_{(t,x) \in (\tau,1] \times \mathbb{R}_+^d} [U(t,x) - V(t,x)] > 0. \tag{2.12}$$

For $\varepsilon \in (0,1)$ let

$$\Phi(t,s,x,y) \doteq U(t,x) - V(s,y) - \frac{1}{\varepsilon}\xi(x-y) - \frac{1}{2\varepsilon}(t-s)^2.$$

Owing to the boundedness of $u$ and $v$ and the non-negativity of $\delta\frac{(1-\tau)^2}{(t-\tau)}$, $\mu$ and $\xi$, $\Phi$ is bounded from above. We also observe that $\Phi$ tends to $-\infty$ as either $s$ or $t$ tends to $\tau$, at a rate that is independent of $(x,y) \in (\mathbb{R}_+^d)^2$. Let the parameter $\sigma \in (0,1)$ satisfy

$$\sigma < \rho \wedge (\delta/4) \wedge (\omega/4M), \tag{2.13}$$

where $M$ is a constant associated with $\xi$ in Lemma 2.1. Since $\sigma > 0$ there exist $(t_\sigma, s_\sigma, x_\sigma, y_\sigma)$ which satisfy

$$\Phi(t_\sigma, s_\sigma, x_\sigma, y_\sigma) \geq \sup_{(t,s,x,y) \in ((\tau,1] \times \mathbb{R}_+^d)^2} \Phi(t,s,x,y) - \sigma \geq \rho - \sigma > 0.$$

We then define the function

$$\Phi_\sigma(t,s,x,y) \doteq \Phi(t,s,x,y) - \frac{\sigma}{2}\left[(t-t_\sigma)^2 + (s-s_\sigma)^2 + \xi(x-x_\sigma) + \xi(y-y_\sigma)\right].$$

Note that if $[(t - t_\sigma)^2 + (s - s_\sigma)^2 + \xi(x - x_\sigma) + \xi(y - y_\sigma)] > 2$ then

$$\begin{aligned}\Phi_\sigma(t, s, x, y) &< \Phi(t, s, x, y) - \sigma \\ &\leq \Phi(t_\sigma, s_\sigma, x_\sigma, y_\sigma) \\ &= \Phi_\sigma(t_\sigma, s_\sigma, x_\sigma, y_\sigma).\end{aligned}$$

We conclude that the supremum in $\sup_{(t,s,x,y) \in ((\tau,1] \times \mathbb{R}_+^d)^2} \Phi_\sigma(t, s, x, y)$ is achieved, and label the point at which it is attained $(\bar{t}, \bar{s}, \bar{x}, \bar{y})$. It follows from $\Phi_\sigma(\bar{t}, \bar{s}, \bar{x}, \bar{y}) > 0$ that

$$\frac{1}{2\varepsilon}(\bar{t} - \bar{s})^2 + \frac{1}{\varepsilon}\xi(\bar{x} - \bar{y}) \leq \sup_{(t,s,x,y) \in ([0,1] \times \mathbb{R}_+^d)^2} (|u(t, x)| + |v(s, y)|) < \infty.$$

Since $\xi(x) \geq m\|x\|^2$, this implies $\bar{t} - \bar{s} \to 0$ and $\bar{x} - \bar{y} \to 0$ as $\varepsilon \to 0$. According to (2.8) and (2.9), $u(t, x) \leq v(t, x)$ for

$$(t, x) \in \Sigma \doteq [[0, 1] \times \{0\}] \cup \left[\{1\} \times \mathbb{R}_+^d\right].$$

We claim that for sufficiently small $\varepsilon > 0$ neither $(\bar{t}, \bar{x})$ nor $(\bar{s}, \bar{y})$ can be in $\Sigma$. Indeed, if this is not true then there exists a subsequence along which both $(\bar{t}, \bar{x})$ and $(\bar{s}, \bar{y})$ converge to the same point in $\Sigma$. Using the continuity of $u$ and $v$, this would violate the strictly positive lower bound $(\rho - \sigma)$ on $\Phi_\sigma$ at $(\bar{t}, \bar{s}, \bar{x}, \bar{y})$, since the limit superior of $\Phi_\sigma(\bar{t}, \bar{s}, \bar{x}, \bar{y})$ along such a subsequence would be bounded above by zero.

We now consider the mapping

$$(t, x) \to u(t, x) - \left[\frac{1}{\varepsilon}\xi(x - \bar{y}) + \frac{1}{2\varepsilon}(t - \bar{s})^2 + \delta\frac{(1 - \tau)^2}{(t - \tau)} + c\mu(x) + \frac{\sigma}{2}(t - t_\sigma)^2 + \frac{\sigma}{2}\xi(x - x_\sigma)\right].$$

Since this function achieves its maximum at $(\bar{t}, \bar{x}) \notin \Sigma$ and $u$ is a subsolution, equation (2.6) must hold. Letting

$$\theta(t, x) \doteq \frac{1}{\varepsilon}\xi(x - \bar{y}) + \frac{1}{2\varepsilon}(t - \bar{s})^2 + \delta\frac{(1 - \tau)^2}{(t - \tau)} + c\mu(x) + \frac{\sigma}{2}(t - t_\sigma)^2 + \frac{\sigma}{2}\xi(x - x_\sigma),$$

we compute

$$\begin{aligned}\theta_t(\bar{t}, \bar{x}) &= \frac{1}{\varepsilon}(\bar{t} - \bar{s}) - \delta\frac{(1 - \tau)^2}{(\bar{t} - \tau)^2} + \sigma(\bar{t} - t_\sigma) \\ D\theta(\bar{t}, \bar{x}) &= \frac{1}{\varepsilon}b(\bar{x} - \bar{y})\Xi(\bar{x} - \bar{y})\bar{n}(\bar{x} - \bar{y}) + \frac{\sigma}{2}b(\bar{x} - x_\sigma)\Xi(\bar{x} - x_\sigma)\bar{n}(\bar{x} - x_\sigma) + cD\mu(\bar{x}).\end{aligned}$$

Now consider $i$ such that $\bar{x}_i = 0$, i.e., $i \in I(\bar{x})$. Since $\bar{y} \in \mathbb{R}_+^d$, it follows that $\langle \bar{x} - \bar{y}, e_i \rangle \leq 0$ for all such $i$, and thus by Condition 2.2 that $\langle \gamma_i, \bar{n}(\bar{x} - \bar{y}) \rangle = \langle \gamma_i, n(z) \rangle \leq 0$ for $i \in I(\bar{x})$ (where $z = \alpha(\bar{x} - \bar{y}), \alpha > 0$, and $z \in \partial B$). The same argument also shows that $\langle \gamma_i, \bar{n}(\bar{x} - x_\sigma) \rangle \leq 0$ for $i \in I(\bar{x})$. Since by Lemma 2.2 $\langle \gamma_i, D\mu(\bar{x}) \rangle < 0$, it follows that $\langle \gamma_i, D\theta(\bar{t}, \bar{x}) \rangle < 0$ for all $i \in I(\bar{x})$.

According to (2.6), it must be true that

$$-\theta_t(\bar{t}, \bar{x}) - H(D\theta(\bar{t}, \bar{x})) \leq 0 \tag{2.14}$$

and $u(\bar{t}, \bar{x}) \leq g(\bar{t}, \bar{x})$. An analogous argument which uses the fact that $v$ is a supersolution shows that either

$$-\zeta_s(\bar{s}, \bar{y}) - H(D\zeta(\bar{s}, \bar{y})) \geq 0, \tag{2.15}$$

with

$$\begin{aligned}
\zeta_s(\bar{s}, \bar{y}) &= \frac{1}{\varepsilon}(\bar{t} - \bar{s}) + \delta\frac{(1-\tau)^2}{(\bar{s}-\tau)^2} - \sigma(\bar{s} - s_\sigma) \\
D\zeta(\bar{s}, \bar{y}) &= \frac{1}{\varepsilon}b(\bar{x} - \bar{y})\Xi(\bar{x} - \bar{y})\bar{n}(\bar{x} - \bar{y}) + \frac{\sigma}{2}b(y_\sigma - \bar{y})\Xi(y_\sigma - \bar{y})\bar{n}(y_\sigma - \bar{y}) - cD\mu(\bar{y})
\end{aligned}$$

or else $v(\bar{s}, \bar{y}) \geq g(\bar{s}, \bar{y})$. Suppose that in fact $v(\bar{s}, \bar{y}) \geq g(\bar{s}, \bar{y})$ occurs infinitely often as $\varepsilon \to 0$. Since $u(\bar{t}, \bar{x}) \leq g(\bar{t}, \bar{x})$, $g$ is continuous, and $\bar{x} - \bar{y} \to 0$ and $\bar{t} - \bar{s} \to 0$, this would imply along the subsequence that

$$\begin{aligned}
\limsup \Phi_\sigma(\bar{t}, \bar{s}, \bar{x}, \bar{y}) &\leq \limsup [u(\bar{t}, \bar{x}) - v(\bar{s}, \bar{y})] \\
&\leq \limsup [g(\bar{t}, \bar{x}) - g(\bar{s}, \bar{y})] \\
&\leq 0.
\end{aligned}$$

Since this contradicts the lower bound $\Phi_\sigma(\bar{t}, \bar{s}, \bar{x}, \bar{y}) \geq \rho - \sigma > 0$, we conclude that (2.14) and (2.15) hold for all sufficiently small $\varepsilon > 0$.

We now subtract (2.14) and (2.15), and use $|\bar{t} - t_\sigma| \vee |\bar{s} - s_\sigma| \leq 2$ and (2.13) to obtain

$$H(D\theta(\bar{t}, \bar{x})) - H(D\zeta(\bar{s}, \bar{y})) \geq -\theta_t(\bar{t}, \bar{x}) + \zeta_s(\bar{s}, \bar{y}) \geq 2\delta - 4\sigma \geq \delta. \tag{2.16}$$

Recall that $L < \infty$ serves as a Lipschitz constant for $u$ and $v$. Hence the characterization of $(\bar{t}, \bar{x})$ and $(\bar{s}, \bar{y})$ as maximizing and minimizing points implies $\|D\theta(\bar{t}, \bar{x})\| \vee \|D\zeta(\bar{s}, \bar{y})\| \leq L$. Let

$$\begin{aligned}
p &\doteq \frac{1}{\varepsilon}b(\bar{x} - \bar{y})\Xi(\bar{x} - \bar{y})\bar{n}(\bar{x} - \bar{y}) \\
v_1 &\doteq \frac{\sigma}{2}b(\bar{x} - x_\sigma)\Xi(\bar{x} - x_\sigma)\bar{n}(\bar{x} - x_\sigma) + cD\mu(\bar{x}) \\
v_2 &\doteq \frac{\sigma}{2}b(y_\sigma - \bar{y})\Xi(y_\sigma - \bar{y})\bar{n}(y_\sigma - \bar{y}) - cD\mu(\bar{y}).
\end{aligned}$$

Recall that $\xi = \Xi^2$, where $\xi$ and $\Xi$ are defined in Lemma 2.1. Then $\xi(\bar{x} - x_\sigma) \leq 2$ implies $\Xi(\bar{x} - x_\sigma) \leq 2$, and therefore $\|\frac{\sigma}{2}b(\bar{x} - x_\sigma)\Xi(\bar{x} - x_\sigma)\bar{n}(\bar{x} - x_\sigma)\| \leq \sigma M$. In a similar fashion one can show $\|\frac{\sigma}{2}b(y_\sigma - \bar{y})\Xi(y_\sigma - \bar{y})\bar{n}(y_\sigma - \bar{y})\| \leq \sigma M$. Also, since $c \leq \omega/4\bar{M}$, where $\bar{M}$ is an upper bound on $\|D\mu(x)\|$, $\|cD\mu(\bar{x})\| \vee \|cD\mu(\bar{y})\| \leq c\bar{M} \leq \omega/4$. Together with (2.13) and the fact that $\omega < 1$, these bounds imply

$$\begin{aligned}
\|v_1\| &\leq \omega \\
\|v_2\| &\leq \omega \\
\|p\| &\leq L + \|v_1\| \leq L + 1.
\end{aligned}$$

Hence (2.11) holds. Since this contradicts (2.16) we conclude that equation (2.12) cannot be true. This implies (2.10), which completes the proof. ∎

# 3  The Differential Game

In this section we formulate a differential game whose value function is a solution to the PDE of the last section. This is done for two reasons. One is to show existence of solutions. The second is so that certain scaling properties of the solution may be verified. In particular, it will show that the solution is linear in radial directions. However, it should be noted that the game we use for these purposes is not the game that motivates our study. The results of this paper will be used elsewhere to construct and analyze optimal service policies in a queueing network, where the optimization problem will be formulated as a game. For any given PDE there are a variety of games whose value function will define a solution to the PDE. We formulate here what seems to be the simplest possible game that will suit our needs. Examples of games that are directly related to queueing network optimization problems can be found in Section 5.

The game will be formulated in the Elliot-Kalton sense [13]. We first give an informal description. In order to properly define the dynamics of the controlled process we must introduce the *Skorokhod Problem*, which is also known in some contexts as the *reflection mapping*. The mapping defined by the solution to the Skorokhod Problem, which will be referred to as the Skorokhod Map, will allow us to construct controlled, constrained processes that are consistent with the boundary conditions used in Section 2. Of course in the queueing problems that motivate the PDE of Section 2 it is the constraints on the state space of the network that define the boundary conditions for the PDE, and not the other way around. We give here the simplest formulation of the SP that will cover our needs. For a more general formulation see [11].

Let
$$C_{\mathbb{R}_+^d}([0,1]:\mathbb{R}^d) \doteq \{\psi \in C([0,1]:\mathbb{R}^d) : \psi(0) \in \mathbb{R}_+^d\},$$
where $C([0,1]:\mathbb{R}^d)$ is the usual space of continuous functions with the sup norm metric. Let a set of vectors that satisfy Condition 2.2 be given. For each point $x$ on the boundary of $\mathbb{R}_+^d$ let
$$d(x) \doteq \left\{ \sum_{i \in I(x)} a_i \gamma_i : a_i \geq 0, \left\| \sum_{i \in I(x)} a_i \gamma_i \right\| = 1 \right\},$$
where as in Section 2 $I(x) \doteq \{i : x_i = 0\}$. The Skorokhod Map assigns to every path $\psi \in C_{\mathbb{R}_+^d}([0,1]:\mathbb{R}^d)$ a path $\phi$ that starts at $\phi(0) = \psi(0)$, but is constrained to $\mathbb{R}_+^d$ as follows. If $\phi$ is in the interior of $\mathbb{R}_+^d$ then the evolution of $\phi$ mimics that of $\psi$, in that the increments of the two functions are the same until $\phi$ hits the boundary of $\mathbb{R}_+^d$. When $\phi$ is on the boundary a constraining "force" is applied to keep $\phi$ in the domain, and this force can only be applied in one of the directions $d(\phi(t))$, and only for $t$ such that $\phi(t)$ is on the boundary. The precise definition is as follows. For $\eta \in C([0,1]:\mathbb{R}^d)$ and $t \in [0,1]$ we let $|\eta|(t)$ denote the total variation of $\eta$ on $[0,t]$ with respect to the Euclidean norm on $\mathbb{R}^d$.

**Definition 3.1** *Let $\psi \in C_{\mathbb{R}_+^d}([0,1]:\mathbb{R}^d)$ be given. Then $(\phi, \eta)$ solves the SP for $\psi$ (with respect to $\mathbb{R}_+^d$ and $\gamma_i, i = 1, ..., d$) if $\phi(0) = \psi(0)$, and if for all $t \in [0,1]$*

1. $\phi(t) = \psi(t) + \eta(t)$,

2. $\phi(t) \in G$,

3. $|\eta|(t) < \infty$,

4. $|\eta|(t) = \int_{[0,t]} 1_{\{\phi(s) \in \partial \mathbb{R}_+^d\}} d|\eta|(s)$,

5. There exists a Borel measurable function $\gamma : [0,1] \to \mathbb{R}_+^d$ such that $d|\eta|$-almost everywhere $\gamma(t) \in d(\phi(t))$, and such that
$$\eta(t) = \int_{[0,t]} \gamma(s) d|\eta|(s).$$

Note that $\eta$ changes only when $\phi$ is on the boundary, and only in the directions $d(\phi)$.

Under Condition 2.2 the Skorokhod Map is Lipschitz continuous on the subset of $C_{\mathbb{R}_+^d}([0,1] : \mathbb{R}^d)$ for which solutions exist [7, 11]. The following condition is sufficient for this subset to be all of $C_{\mathbb{R}_+^d}([0,1] : \mathbb{R}^d)$. This condition, which is rather weak, is known in the literature as the *completely-$\mathcal{S}$* condition ([19]). Within the setup of the present paper this condition is also necessary for the existence of solutions ([4], [11, Section 4 of part I]). Note that Condition 3.1 implies Condition 2.3 of Section 2.

**Condition 3.1** *Given $\kappa \subset \{1, ..., d\}$ there exist $b_i \geq 0, i \in \kappa$ such that if $\gamma \doteq \sum_{i \in \kappa} b_i \gamma_i$, then $\langle \gamma, e_i \rangle > 0$ for all $i \in \kappa$.*

We next define a constrained ordinary differential equation. As is proved under Conditions 2.2 and 3.1 in [7], one can define a projection $\pi : \mathbb{R}^d \to \mathbb{R}_+^d$ that is consistent with the constraint directions $\{\gamma_i, i = 1, ..., d\}$, in that $\pi(x) = x$ if $x \in \mathbb{R}_+^d$, and if $x \notin \mathbb{R}_+^d$ then $\pi(x) - x = \alpha r$, where $\alpha \geq 0$ and $r \in d(\pi(x))$. Figure 1 illustrates the projection for a two dimensional problem.

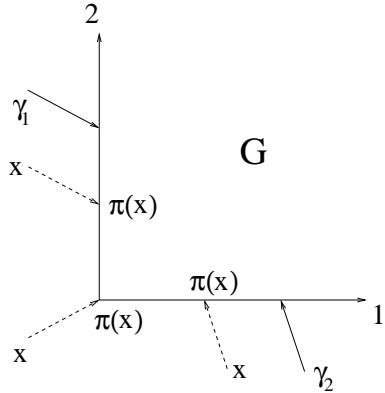

Figure 1: The discrete projection

With this projection given, we can now define for each point $x \in \partial \mathbb{R}_+^d$ and each $v \in \mathbb{R}^d$ the *projected velocity* ([5])
$$\pi(x,v) \doteq \lim_{\Delta \downarrow 0} \frac{\pi(x + \Delta v) - \pi(x)}{\Delta}.$$

For details on why this limit is always well defined and further properties of the projected velocity we refer to [5, Section 3 and Lemma 3.8]. Let $\beta : [0,1] \to \mathbb{R}^d$ have the property that each component of $\beta$ is integrable over each interval $[0,T], T < \infty$. Then the ODE of interest is given by

$$\dot\phi(t) = \pi(\phi(t), \beta(t)). \tag{3.1}$$

An absolutely continuous function $\phi : [0,1] \to \mathbb{R}^d_+$ is a solution to (3.1) if the equation is satisfied in an a.e. sense in $t$. By using the regularity properties of the associated Skorokhod Map, one can prove that all the standard qualitative properties (existence and uniqueness of solutions, stability with respect to perturbations, etc.) hold [7, 10]. In fact, because of the particularly simple nature of the right hand side (i.e., $\pi(\phi(t), \beta(t))$ rather than $\pi(\phi(t), b(\phi(t)) + \beta(t))$ for some function $b$), one can show that $\phi$ solves (3.1) if and only if $\phi$ is the image of $\psi(t) \doteq \int_0^t \beta(s)ds + x$ under the Skorokhod Map, in which case all such issues become trivial.

The ODE (3.1) will define the dynamics in the game we consider. In order to complete this formal description of the game, we must also describe a cost structure. Associated with each of the convex functions $H_j, j = 1, ..., J$, is its Legendre transform, which is defined for each $\beta \in \mathbb{R}^d$ by

$$L_j(\beta) \doteq \sup_{p \in \mathbb{R}^d} \left[ \langle \beta, p \rangle - H_j(p) \right].$$

As is well known ([20]), each function $L_j : \mathbb{R}^d \to \mathbb{R} \cup \{\infty\}$ is never equal to $-\infty$, convex on the set $\{\beta : L_j(\beta) < \infty\}$, and lower semicontinuous. The maximizing player in the game will select a measurable control $\beta(t)$, which will control the dynamics through (3.1). The minimizing player will select a measurable control $k(t)$, which will take values in the set $\{1, ..., J\}$. This control will not directly affect the dynamics, but instead chooses which of the costs $L_k(\beta)$ will be applied to the maximizing player. If ever $\phi(t) = 0$ then play is stopped and a stopping cost of zero is paid. In addition, the minimizing player will choose a control $\lambda(t)$ that takes values in $\{0, 1\}$. If $t$ is the first time that $\lambda(t) = 1$, then play is stopped and a stopping cost of $g(t, \phi(t))$ is incurred. If at time 1 play has not been stopped either by the minimizing player or by hitting the origin, then a terminal cost of $f(\phi(1))$ is paid. Hence the cost function for controls $(\beta, k, \lambda)$ and initial condition $\phi(t) = x$ at time $t$ is defined as

$$C_t(x, \beta, k, \lambda) \doteq \left[ -\int_t^{\tau \wedge \sigma \wedge 1} L_{k(s)}(\beta(s))ds + 1_{\{\tau < \sigma, \tau < 1\}} g(\tau, \phi(\tau)) + 1_{\{\tau \geq 1, \sigma > 1\}} f(\phi(1)) \right],$$

where $\sigma \doteq \inf\{s : \phi(s) = 0\}$ and $\tau \doteq \inf\{s : \lambda(s) = 1\}$.

Although this completes the description of the dynamics and the cost structure, it does not fully describe the game, since we have not indicated which player has the "information advantage," nor have we indicated how that will be modeled. For the game we are interested in it is the maximizing player that is given the advantage [13]. Hence we consider the corresponding *upper game*. (Interestingly, for the games that are directly related to queueing optimization problems the value for the upper and lower games often coincide, and both equal the value of this upper game.)

Let

$$\mathcal{J}_t \doteq \{k : [t,1] \to \{1, ..., J\} : k \text{ is Borel measurable}\}$$

and

$$\mathcal{B}_t \doteq \{\beta : [t,1] \to \mathbb{R}^d : \beta \text{ is Borel measurable}\}.$$

A *strategy for the maximizing player at time t* is a map $\nu : \mathcal{J}_t \to \mathcal{B}_t$ with the property that $\nu$ is *nonanticipating*: if $k_1(t) = k_2(t)$ for a.e. $t \in [s, r]$, then $\nu[k_1](t) = \nu[k_2](t)$ for a.e. $t \in [s, r]$. Let $\Lambda_t$ be the set of mappings from $[t, 1]$ to $\{0, 1\}$. Owing to the nonanticipating nature of strategies it would make no difference if we allowed the strategy for the maximizing player to also depend on the minimizing player's stopping control, and so to simplify this dependence is omitted. Let $\Gamma_t$ denote the collection of all strategies for the maximizing player at time $t$. Then the value function for the game is defined by

$$V(x,t) \doteq \sup_{\nu \in \Gamma_t} \inf_{k \in \mathcal{J}_t, \lambda \in \Lambda_t} C_t(x, \nu[k], k, \lambda). \tag{3.2}$$

In this paper we consider only Lipschitz continuous solutions. We therefore assume a compatibility condition between the stopping cost and the terminal cost, which is that $f(x) = g(1, x)$ for $x \in \mathbb{R}_+^d$. This condition will be satisfied by the time dependent problems used to analyze the time independent problem in Section 4.

**Theorem 3.2** *Assume Conditions 2.1, 2.2, and 3.1, and that $f(x) = g(1, x)$. Let $V(x, t)$ be defined by (3.2), and assume that $V$ is Lipschitz continuous. Then $V(x, t)$ is the unique viscosity solution to (2.2)–(2.5).*

As is the usual case, the proof of this theorem follows from a dynamic programming principle. The following theorem can be proved using standard arguments [14].

**Theorem 3.3** *Assume Conditions 2.1, 2.2 and 3.1, and that $f(x) = g(1, x)$. Let $V(t, x)$ be defined by (3.2), and assume that $V$ is Lipschitz continuous. Let $0 \leq t < 1$ and consider any point $x \in \mathbb{R}_+^d \setminus \{0\}$ for which $V(t, x) < g(t, x)$. Then there exists $\delta_0 \in (0, 1 - t)$ such that for any $\delta \in (0, \delta_0)$*

$$V(t, x) = \sup_{\nu \in \Gamma_t} \inf_{k \in \mathcal{J}_t} \left[ -\int_t^{t+\delta} L_{k(s)}(\beta(s)) ds + V(t + \delta, \phi(t + \delta)) \right],$$

*where $\phi$ solves (3.1) with $\phi(t) = x$.*

**Proof of Theorem 3.2.** Properties (2.8) and (2.9) follow directly from the definition of $V(t, x)$. Thus we need only consider conditions (2.6) and (2.7).

We first show that (2.7) holds. Let $\theta : (0, T) \times \mathbb{R}_+^d \to \mathbb{R}$ be a continuously differentiable function, and let $(s, y) \in (0, T) \times \left( \mathbb{R}_+^d \setminus \{0\} \right)$ be a local minimum of $V - \theta$. We can assume without loss that $V(s, y) = \theta(s, y)$. If $V(s, y) \geq g(s, y)$ then (2.7) is satisfied, and so we must only consider the case when $V(s, y) < g(s, y)$. If (2.7) is not true, then there must be $a > 0$ such that

$$\theta_t(s, y) + H(D\theta(s, y)) \geq a \tag{3.3}$$

and

$$\langle D\theta(s, y), \gamma_i \rangle \geq a \tag{3.4}$$

for all $i \in I(y)$. Let $\delta_0$ satisfy the conclusion of Theorem 3.3 for $(t,x) = (s,y)$.

The definition of $H$ and (3.3) imply

$$\theta_t(s,y) + \wedge_{j=1}^{J} \sup_{\beta \in \mathbb{R}^d} [\langle \beta, D\theta(s,y) \rangle - L_j(\beta)] \geq a.$$

For $j \in \{1, ..., J\}$ choose $\bar{\beta}_j \in \mathbb{R}^d$ so that

$$\theta_t(s,y) + [\langle \bar{\beta}_j, D\theta(s,y) \rangle - L_j(\bar{\beta}_j)] \geq \frac{a}{2}. \tag{3.5}$$

These vectors are used to define a strategy (in $\Gamma_s$) for the maximizing player by setting

$$\nu[k](t) = \bar{\beta}_{k(t)}$$

for any control $k \in \mathcal{J}_s$. Let $\phi$ denote the controlled dynamics when this strategy is used and the initial condition is $\phi(s) = y$. Since the mapping $y \to I(y)$ is upper semicontinuous (i.e., given any $z \in \mathbb{R}^d$, there is $\varepsilon > 0$ such that $I(y) \subset I(z)$ for all $y \in \mathbb{R}^d$ satisfying $\|y - z\| \leq \varepsilon$), we can assume that $I(\phi(r)) \subset I(y)$ for all $r \in [0, \delta_0]$. Suppose for a given $r \in [0, \delta_0]$ that $k(r) = j$. In this case $\bar{\beta}_j$ can be related to $\phi(r)$ in a variety of ways, depending on the value of $I(\phi(r))$. However, since in all cases $I(\phi(r)) \subset I(y)$, the definition of the Skorokhod Problem implies that

$$\dot{\phi}(r) = \bar{\beta}_j + \sum_{i \in I(y)} a_i \gamma_i$$

for some numbers $a_i \geq 0, i \in I(y)$. Hence by (3.4),

$$\theta_t(s,y) + \langle D\theta(s,y), \bar{\beta}_j \rangle + \sum_{i \in I(y)} \langle D\theta(s,y), a_i \gamma_i \rangle \geq \theta_t(s,y) + \langle D\theta(s,y), \bar{\beta}_j \rangle.$$

When combined with (3.5) and the smoothness of $\theta$, the last display implies

$$\frac{d\theta(r, \phi(r))}{dr} - L_{k(r)}(\bar{\beta}_{k(r)}) \geq \frac{a}{4}$$

for almost all $r \in [0, \delta]$ if $\delta \in (0, \delta_0)$ is sufficiently small. Hence for sufficiently small $\delta \in (0, \delta_0)$ and any control $k$

$$\theta(s+\delta, \phi(s+\delta)) - \theta(s,y) - \int_s^{s+\delta} L_{k(r)}(\nu[k](r)) dr \geq \frac{a\delta}{4} \tag{3.6}$$

However, by Theorem 3.3

$$V(s,y) \geq \inf_{k \in \mathcal{J}_s} \left[ -\int_s^{s+\delta} L_{k(r)}(\beta(r)) dr + V(s+\delta, \phi(s+\delta)) \right].$$

Hence there is a control $k \in \mathcal{J}_s$ such that

$$V(s,y) \geq -\int_s^{s+\delta} L_{k(r)}(\beta(r)) dr + V(s+\delta, \phi(s+\delta)) - \frac{a\delta}{8}.$$

Subtracting the last display from (3.6) and using $V(s, y) = \theta(s, y)$ gives

$$V(s + \delta, \phi(s + \delta)) - \theta(s + \delta, \phi(s + \delta)) \leq -\frac{a\delta}{8}$$

for all sufficiently small $\delta \in (0, \delta_0)$. This contradicts the fact that $(s, y)$ is a local minimum of $V - \theta$, and so it follows that (2.7) is valid.

To complete the proof we must show that (2.6) holds. Let $\theta : (0, T) \times \mathbb{R}^d_+ \to \mathbb{R}$ be continuously differentiable, and let $(s, y) \in (0, T) \times \left(\mathbb{R}^d_+ \backslash \{0\}\right)$ be a local maximum of $V - \theta$. Once again we can assume that $V(s, y) = \theta(s, y)$. Since the minimizing player can always stop the game it follows that

$$V(s, y) - g(s, y) \leq 0.$$

Thus if (2.6) is not true then there is $a > 0$ such that

$$\theta_t(s, y) + H(D\theta(s, y)) \leq -a \tag{3.7}$$

and

$$\langle D\theta(s, y), \gamma_i \rangle \leq -a \tag{3.8}$$

for all $i \in I(y)$. Let $\delta_0$ satisfy the conclusion of Theorem 3.3 for $(t, x) = (s, y)$.

Let $\bar{j} \in \{1, ..., J\}$ be a minimizer in $\wedge_{j=1}^J H_j(D\theta(s, y))$. We define a control for the minimizing player by setting $k(r) = \bar{j}$ for all $r \in [0, \delta_0]$. Since this control is constant, a strategy for the maximizing player is just an open loop control $\beta(r)$ (which can depend on $x$ and $\bar{j}$). According to (3.7),

$$\theta_t(s, y) + \langle D\theta(s, y), \beta(r) \rangle - L_{k(r)}(\beta(r)) \leq -a \tag{3.9}$$

for almost all $r \in [0, \delta_0]$.

Let $\varepsilon > 0$ be such that $\|y - z\| < \varepsilon$ implies $I(z) \subset I(y)$. First assume that $\|\phi(r) - y\| \geq \varepsilon$ for some $r \in [0, \delta]$. Since $H_{\bar{j}}(\alpha) < \infty$ for all $\alpha \in \mathbb{R}^d$, $L_{\bar{j}}$ grows superlinearly:

$$\liminf_{c \to \infty} \inf_{\beta : \|\beta\| = c} L_{\bar{j}}(\beta)/c = \infty.$$

Also, there is $M < \infty$ such that for all $x \in \mathbb{R}^d_+$ and $v \in \mathbb{R}^d$

$$\|\pi(x, v)\| \leq M(1 + \|v\|).$$

Elementary estimates using the last two displays show that if $\delta \in (0, \delta_0)$ is sufficiently small then whenever $\phi$ travels more than $\varepsilon$ away from $y$ in time $\delta$,

$$\theta(s + \delta, \phi(s + \delta)) - \theta(s, y) - \int_s^{s+\delta} L_{k(r)}(\nu[k](r))dr \leq -\frac{a\delta}{4} \tag{3.10}$$

Next consider the case $\|\phi(r) - y\| < \varepsilon$ for all $r \in [0, \delta]$, where $\delta \in (0, \delta_0)$. Then for almost every such $r$ there exist numbers $a_i \geq 0, i \in I(y)$ such that

$$\dot{\phi}(r) = \beta(r) + \sum_{i \in I(y)} a_i \gamma_i.$$

By choosing $\delta \in (0, \delta_0)$ smaller if need be, the smoothness of $\theta$, (3.9) and (3.8) imply

$$\begin{aligned}\frac{d\theta(r, \phi(r))}{dr} - L_{k(r)}(\bar{\beta}_{k(r)}) &= \theta_t(r, \phi(r)) + \langle D\theta(r, \phi(r)), \bar{\beta}_j \rangle \\ &\quad + \sum_{i \in I(y)} \langle D\theta(r, \phi(r)), a_i \gamma_i \rangle - L_{k(r)}(\bar{\beta}_{k(r)}) \\ &\leq -a/2.\end{aligned}$$

Hence for sufficiently small $\delta \in (0, \delta_0)$ and all $\beta(\cdot)$ (3.10) holds. Arguing as before, this estimate and Theorem 3.3 contradict the fact that $(s, y)$ is a local maximum of $V - \theta$. ∎

## 4 A Simple Characterization of the Value of an Infinite Time Game

In the previous sections we proved a uniqueness result and representation theorem for the PDE associated to a finite time game. In this section we consider an infinite time analogue of this game. The value function for this game is of prime interest, in that it will play a key role in characterizing optimal robust service policies for certain network problems. It is analogous to the *cost potential* or *storage function* in standard deterministic robust control.

One can formally write down a PDE that should be satisfied by value functions of this sort simply by deleting the $V_t$ term in the PDE of Section 2 and assuming that $g$ is independent of $t$. However, it turns out (as is also the case in the usual robust control setting, see [17]) that such equations have poor uniqueness properties, even in the viscosity sense. Hence we take a different tack. Let $V(x)$ denote the value for the infinite time problem, and let $V(t, x)$ denote the value defined in the last section. We will use the fact that if $t < 1$ is fixed, then $V(x)$ and $V(t, x)$ *coincide* as functions of $x$ for all $x$ in a neighborhood of the origin. Under conditions stronger than those of the last two sections, we will prove (for an appropriate choice of the stopping cost $g$) that $V(x)$ is entirely characterized by $V(t, x)$, and, in particular, by $V(0, x)$. This will allow the fact that $V(t, x)$ is the unique solution to a PDE to be brought into consideration, and ultimately lead to a simple analytic characterization of $V(x)$.

We begin by defining the infinite time game. As in Section 2, the dynamics of the game are defined by the constrained ODE (3.1), except that solutions are now considered on the infinite time interval $[0, \infty)$. We wish to define a game where play is not stopped until the state of the process hits the origin. Even though there is no fixed time by which the origin should be hit, we consider only controls for the minimizing player that will force this to happen. In particular, and in contrast with the game of Section 3, there is no other way for the minimizing player to stop the game.

The game which satisfies these requirements can be defined as follows. Let

$$\mathcal{J} \doteq \{k : [0, \infty) \to \{1, ..., J\} : k \text{ is Borel measurable}\}$$

and

$$\mathcal{B} \doteq \{\beta : [0, \infty) \to \mathbb{R}^d : \beta \text{ is Borel measurable}\},$$

and let $\Gamma$ denotes the set of non-anticipating maps from $\mathcal{J}$ to $\mathcal{B}$. We then define

$$V(x) \doteq \sup_{\nu \in \Gamma} \inf_{k \in \mathcal{J}: \sigma < \infty} C_0^\sigma(x, \nu[k], k), \tag{4.1}$$

where $\phi$ is defined by (3.1),

$$C_0^\sigma(x, \beta, k) \doteq \left[ -\int_0^\sigma L_{k(t)}(\beta(t)) dt \right],$$

$\sigma \doteq \inf\{t : \phi(t) = 0\}$, and where the infimum over the empty set is $\infty$.

One can easily verify by a time rescaling argument that $V$ is *linear in radial directions*, which means that for all $x \in \mathbb{R}_+^d$ and all $\alpha > 0$

$$V(\alpha x) = \alpha V(x).$$

In this section we use the following additional condition.

**Condition 4.1** *Each function $L_j, j \in \{1, ..., J\}$ takes the form $\hat{L}_j(\beta) - 1$, where*

$$\hat{L}_j(\beta) \doteq \begin{cases} 0 & \text{if } x \in C_j \\ \infty & \text{otherwise,} \end{cases}$$

*with each $C_j$ a closed nonempty convex set.*

This condition is satisfied by Examples 1 and 2 of Section 5, but not by Example 3. At the present time it is not known if one can generalize Theorem 4.3 to cover Example 3 as well. The game that defines $V(x)$ has a particularly simple interpretation when Condition 4.1 holds. The cost $\hat{L}_j(\beta) = \infty$ for $\beta \notin C_j$ acts as a constraint on the maximizing player, forcing $\beta(t) \in C_j$ whenever the minimizing player chooses $k(t) = j$. Given that this constraint is satisfied a.s. in $t$, the cost $C_0^\sigma(x, \beta, k)$ is simply the time that it takes to hit the origin, given a starting position $x$. Hence the goal of the maximizing player is to delay this event for as long as possible. The minimizing player wishes for this time to be as small as possible, and attempts to achieve this by selecting the constraint on $\beta(t)$ appropriately.

We have observed that $V(x)$ represents the optimal time till the origin is reached. Let $g(x,t) \doteq (1-t)$ and $R \doteq \{x : V(x) \leq 1\}$. We claim that when the stopping cost $g$ is used in the game of Section 3 the solution is given for $t \in [0,1]$ by

$$V(t, x) = \begin{cases} V(x) & \text{if } x \in (1-t)R, \\ (1-t) & \text{otherwise.} \end{cases} \tag{4.2}$$

The precise statement is as follows.

**Theorem 4.1** *Let $g(t, x) \doteq (1 - t)$ and $f(x) \doteq 0$. Assume that $V$ defined by (4.1) is finite and Lipschitz continuous, and that Conditions 2.1, 2.2, 3.1 and 4.1 are satisfied. Then $V(t, x)$ defined by (4.2) is the unique viscosity solution to (2.2)–(2.5).*

**Proof** We give the proof for $t = 0$. The general case uses the same argument. Fix $x \in R^\circ$ and let $\varepsilon > 0$ be given. In the definition of $V(x)$, choose $\nu$ such that

$$\inf_{k \in \mathcal{J}: \sigma < \infty} C_0^\sigma(x, \nu[k], k) \geq V(x) - \varepsilon.$$

This means for all $k \in \mathcal{J}$ that $\sigma \geq V(x) - \varepsilon$. We adapt $\nu$ in the obvious way for use in the definition of the finite time problem $V(0, x)$. Because $\nu$ is non-anticipating this adaptation is well defined. If the minimizing player in $V(0, x)$ never uses the stopping control, then they pay at least $V(x) - \varepsilon$ if $\sigma \leq 1$ and $1$ if $\sigma > 1$. If they do stop play at some time $s \in [0, 1]$, then the cost will be at least $s + g(s, \phi(s)) = s + (1 - s) = 1$. Since $x \in R$, the radial linearity of $V(x)$ and the definition of $R$ imply $1 \geq V(x)$, and so in all cases the cost is at least $V(x) - \varepsilon$. Since $\varepsilon > 0$ is arbitrary, $V(0, x) \geq V(x)$.

Next let $\nu$ be any strategy that appears in the definition of $V(0, x)$. Given $\varepsilon > 0$ we must find a control $(k, \lambda)$ for the minimizing player such that

$$C_0(x, \nu[k], k, \lambda) \leq V(x) + \varepsilon. \tag{4.3}$$

Since $x \in R^\circ$ we can find $\varepsilon > 0$ such that $V(x) + \varepsilon < 1$. Let $\bar{\nu}$ be any non-anticipating extension of $\nu$ to trajectories $k \in \mathcal{J}$ that are defined on the interval $[0, \infty)$. Choose $k \in \mathcal{J}$ in the definition of $V(x)$ such that $\sigma \leq V(x) + \varepsilon$, and by abusing notation let $k$ also denote the corresponding restriction to $[0, 1]$. Note that $\sigma \leq V(x) + \varepsilon$ implies $\sigma < 1$. By choosing to never exercise the stopping control, the minimizing player obtains (4.3). Since $\varepsilon > 0$ is arbitrary we find that $V(0, x) \leq V(x)$.

Finally we must consider $x \notin R^\circ$. Since the minimizing player can always stop play immediately, $V(0, x) \leq 1$. On the other hand, the same argument as in the case $x \in R^\circ$ shows that $V(0, x) \geq 1$, and so $V(0, x) = 1$. ∎

The definition of viscosity solutions in section 2 is phrased in terms of test functions. An alternative definition in terms of sub- and superdifferentials is standard (see [3, 15]). The following result from [3] provides a link between the two definitions. It is stated in [3, Lemma II.1.7] for an open domain, but the proof of the following version is identical.

**Lemma 4.2** *Let $\Omega$ denote $(0, 1) \times \left(\mathbb{R}_+^d \setminus \{0\}\right)$ and let $w \in C(\Omega)$. Then $p \in D^+ w(x)$ [resp., $p \in D^- w(x)$] if and only if there exists $\theta \in C^1(\Omega)$ such that $D\theta(x) = p$ and $w - \theta$ has a local maximum [resp., minimum] at $x$.*

As a result, the following statement is equivalent to the definition of Section 2. A bounded Lipschitz continuous function $w : [0, 1] \times \mathbb{R}_+^d \to \mathbb{R}$ is a viscosity subsolution [resp., supersolution] to (2.2)–(2.5) if the following conditions hold. If $(a, p) \in D^+ V(s, y), s \in (0, 1), y \in \mathbb{R}_+^d \setminus \{0\}$ then

$$\max\left\{\left([-a - H(p)] \wedge \min_{i \in I(y)}[-\langle p, \gamma_i \rangle]\right); w(s, y) - g(s, y)\right\} \leq 0 \tag{4.4}$$

[resp., if $(a, p) \in D^- V(s, y), s \in (0, 1), y \in \mathbb{R}_+^d \setminus \{0\}$ then

$$\max\left\{\left([-a - H(p)] \vee \max_{i \in I(y)}[-\langle p, \gamma_i \rangle]\right); w(s, y) - g(s, y)\right\} \geq 0]. \tag{4.5}$$

Moreover, (2.8) and (2.9) hold.

The following analytic characterization of $V$ follows from Theorem 4.1. Note that even when $V$ is smooth the set of sub and superdifferentials is multivalued on $\partial \mathbb{R}_+^d$. In particular, if $p \in D^+V(x)$ [resp., $p \in D^-V(x)$] and $x \in \partial \mathbb{R}_+^d$, then $p + \sum_{i \in I(x)} a_i e_i \in D^+V(x)$ [resp., $p - \sum_{i \in I(x)} a_i e_i \in D^-V(x)$] for all $a_i \geq 0, i \in I(x)$.

**Theorem 4.3** *Assume that $V$ defined by (4.1) is finite and Lipschitz continuous, and that Conditions 2.1, 2.2, 3.1 and 4.1 are satisfied. Then $V$ satisfies the following conditions:*

$$p \in D^-V(x), x \in R_+^d \setminus \{0\} \Rightarrow H(p) \wedge \min_{i \in I(x)} \langle p, \gamma_i \rangle \leq 0, \tag{4.6}$$

$$V(0) = 0, \tag{4.7}$$

$$p \in D^+V(x), x \in R_+^d \setminus \{0\}, \text{ and } b \in [0,1] \Rightarrow [H(bp) - (1-b)] \vee \max_{i \in I(x)} \langle bp, \gamma_i \rangle \geq 0, \tag{4.8}$$

*Conversely, suppose that a Lipschitz continuous and radially linear function $V : \mathbb{R}_+^d \to \mathbb{R}$ is given that satisfies conditions (4.6)–(4.8), and such that the level set $R \doteq \{x : V(x) \leq 1\}$ is bounded. Then $V$ equals the value function defined by (4.1).*

**Proof** Assume that the function $V(x)$ defined by (4.1) is finite and Lipschitz continuous. Property (4.7) follows directly from the definition of $V$. Define $V(t,x)$ by (4.2). Equations (4.6) and (4.8) will be proved by applying Theorems 3.2 and 4.1, and then using the alternative phrasing of conditions (2.6)–(2.9) in terms of sub and superdifferentials.

First assume that $x \in [(1-t)R]^\circ$, where $t \in (0,1)$. For such points the definition of $R$ implies that $V(t,x) < g(t,x)$. If $p \in D^-V(x)$ and $x \in \mathbb{R}_+^d \setminus \{0\}$, then by (4.2) $(0,p) \in D^-V(t,x)$. Inserting this into equation (4.5) and using $V(t,x) < g(t,x)$ gives (4.6). Next let $x \in \partial[(1-t)R]$ and $t \in (0,1)$ be given. Then

$$D^+V(t,x) = \{(-(1-b), bp) : p \in D^+V(x), b \in [0,1]\}.$$

Together with (4.4), the last display gives (4.8).

To prove the last sentence of the theorem, let $V(x)$ be a Lipschitz continuous function and radially linear function that satisfies conditions (4.6)–(4.8). Define $V(t,x)$ by (4.2). We will verify the sub and superdifferential forms of (2.6)–(2.9). Conditions (2.8) and (2.9) follow directly from (4.7), the fact that $f(x) = 0$, and the definition of $V(t,x)$ via (4.2). Next suppose that $(a,p) \in D^+V(s,y)$, with $s \in (0,1)$ and $y \in \mathbb{R}_+^d \setminus \{0\}$. If $y \in ((1-s)R)^\circ$ then the definition of $V(t,x)$ implies $V(s,y) < g(s,y)$, $a = 0$, and $p \in D^+V(y)$. In this case (4.8) with $b = 1$ implies (4.4). If $y \in \partial((1-s)R)$ then $V(s,y) = g(s,y)$ and $(a,p) = (-(1-b), bq)$ for some $q \in D^+V(y)$, and again (4.8) implies (4.4). If $y \notin (1-s)R$ then $V(s,y) = g(s,y)$ and $(a,p) = (-1,0)$, and (4.4) follows from the fact that under Condition 4.1

$$H(0) = \wedge_{j=1}^J \sup_{\beta \in \mathbb{R}^d} -L_j(\beta) = 1.$$

Lastly we turn to condition (4.5). Suppose that $(a,p) \in D^-V(s,y)$, with $s \in (0,1)$ and $y \in \mathbb{R}^d_+ \setminus \{0\}$. If $y \in ((1-s)R)^\circ$ then $V(s,y) < g(s,y)$, and so $a = 0$ and $p \in D^-V(y)$. In this case (4.6) implies (4.5). If $y \notin (1-s)R$ then $V(s,y) = g(s,y)$, and so (4.5) holds as well. Having verified all of conditions (2.6)–(2.9), the proof of the theorem is complete. ∎

## 5 Three Examples From Queueing Network Optimization

In this section we present three examples. The aim is to show how the game defined in Section 3 can be used to formulate and study problems in the optimal control of queueing networks. The examples all consider the same 4 class, 2 station queueing network, but with different cost criteria, motivations, and interpretations. Although we focus here on service control, the ideas can be applied to problems with routing control as well.

We first describe the dynamical model that will be used for all three examples. The queueing system consists of four queues, labeled 1, 2, 3 and 4, and two servers, labeled A and B (see Figure 2). Customers enter queue 1 at rate $\lambda \in [0,\infty)$. Queues 1 and 4 (resp., 2 and 3) are served by server A (resp., B). Queue $i$ is served at rate $\mu_i$ when the respective server serves only this queue. The control in this problem is the choice of which class is served, and hence the control takes 4 values. It will turn out that this control problem has the same value as the corresponding "relaxed" version, which allows that server to split service effort between customer classes.

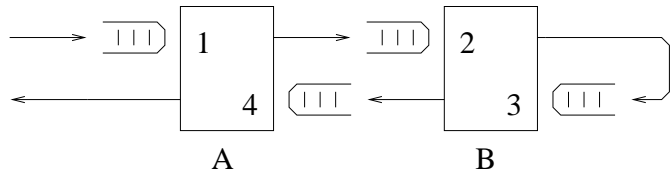

Figure 2: A system with four queues and two servers

For particular fixed choices of the control policy, this system is studied for a deterministic model in [16] and for a stochastic model in [18]. It has been shown that this system can experience instability for non-idling service policies, even if the usual load conditions are satisfied. This nontrivial behavior with regard to stability criteria makes the example interesting from the point of view of control and robust stabilization.

We next describe the Skorokhod Problem model that is appropriate for this system. Our "derivation" of the Skorokhod Problem is intended only as a heuristic. A more careful analysis can be applied to relate a variety of more detailed models to the model we use via a law of large numbers type approximation. At the level of modeling we consider (which is sometimes called the "fluid" level), the system is described by a continuous time, continuous state dynamical model. The state $\phi \in \mathbb{R}^4_+$ represents the queue lengths at sites 1, 2, 3 and 4. The role of the Skorokhod Problem is to apply the proper compensation to the state when a service is attempted at an empty queue. A key observation here is that this is the *only* role the Skorokhod Problem will play. In particular, the Skorokhod Problem will not be required to enforce non-idling behavior. In the setting of multiclass problems, it is common when considering issues of performance analysis to restrict to a fixed non-idling service policy, which then must be reflected in the form of the Skorokhod Problem [12]. Unfortunately, this can lead to a Skorokhod Problem that is not well defined, in that there can be multiple solutions for a single input $\psi$ [18]. In the present paper we allow the service to be a control, which turns out to be extremely convenient in that it *always* leads to a well defined Skorokhod Problem. We will return to this point after the description of the Skorokhod Problem is complete.

Since the state space has been identified, all that remains in defining the Skorokhod Problem is to identify the vectors $\gamma_i$. To see how this can be accomplished, we first consider the evolution of the system away from any of the boundaries, i.e., at points where all queues are nonempty. The dynamics will of course depend on the current control. Let $\mathcal{U} \doteq \{(1,2),(1,3),(4,2),(4,3)\}$ denote the control space. Hence if the control $u = (s_A, s_B) \in \mathcal{U}$ is chosen, the server $A$ serves $s_A$ and server $B$ serves $s_B$. Suppose, for example, that $u = (1,3)$. Then the dynamics away from the boundary with this choice of control are

$$\dot{\phi} = \lambda e_1 + \mu_1(e_2 - e_1) + \mu_3(e_4 - e_3) = (\lambda - \mu_1, \mu_1 - \mu_2, \mu_2, 0).$$

The interior dynamics with the other choices of the control are obvious analogues.

We next consider the behavior when one of the queues is empty. For $i = 1, 2, 3, 4$, the constraint mechanism on the face $\{x \in \mathbb{R}_+^4 : I(x) = \{i\}\}$ must compensate for any service that is attempted when the queue is empty. As noted previously we will *not* assume, a priori, that the controls are restricted to non-idling service policies. In fact, the controls that appear in the definition of the value functions are "open loop," which simply means that they are functions of time but not of the state of the system (see Section 3). It is certainly possible that some of these controls will allow parts of the system to idle even when there is work to be done. We will rely on optimality to indirectly eliminate such controls when they are not appropriate. Hence if control $u = (1,3)$ is selected, and if queues 1 and 3 are empty, the servers continue to (attempt to) serve queues 1 and 3, and ignore queues 2 and 4. Suppose that $\phi_1$ is zero, and that all other queues are nonempty. Then regardless of the current state of the control, if a service were attempted to queue 1 then the compensating action is to increase queue 1 in proportion to the erroneous service, and to also reduce the length of queue 2 by the same amount. Thus the direction of constraint on this face is $\gamma_1 = e_1 - e_2$. A similar argument applies at the faces which correspond to queues 2 and 3 being empty, and so $\gamma_2 = e_2 - e_3$ and $\gamma_3 = e_3 - e_4$. Since the output of queue 4 is not fed into any other queue, $\gamma_4 = e_4$.

Those familiar with the Skorokhod Problem will recognize that we have obtained the same Skorokhod Problem as the one used to model a single class tandem queueing network with four stations. This means that by optimizing over the service policy and by using open loop controls, we effectively decouple the different customer classes, at least from the point of view of the Skorokhod Problem and related analytic machinery (e.g., the PDE of Sections 2 and 4). This significantly simplifies the analysis, and to the authors' knowledge this fact has not been observed previously. The situation is very much analogous to the classical theory of controlled unconstrained ODEs. Here the system might be modeled by the ODE $\dot{\phi} = b(\phi, u)$, with $b(\cdot, \cdot)$ a smooth function of its arguments and $u$ the control. If one considers a fixed feedback control then $u$ may not be smooth, and there can be difficulty in associating a well defined dynamical model. If one uses open loop controls then there is no difficulty at all in defining quantities such as the system dynamics and the value function. Of course if one is primarily interested in feedback controls then in a sense the difficulty is only deferred, since one must eventually construct a feedback control that performs in at least a nearly optimal fashion. Nonetheless, the use of open loop controls (or open loop controls and strategies in the game case) is important, since quantities such as the value function can be properly defined and used to study the problem further.

With the directions of constraint now identified, the system dynamics are now defined as in Section 3 through equation (3.1). Conditions 2.2 and 3.1 are easily verified for this model. See, for example, [7, 11]. All that remains is to specify the cost structure, and it is here that our examples differ. We begin with the simplest example.

**Example 1: Optimal Stabilizing Service Policy.** We wish to define a problem in which the arrival and service rates have fixed values, and the issue is to design a control strategy such that for any given nonzero initial state, the process $\phi$ is returned to the origin as quickly as possible. The index set $\{1, ..., J\}$ of the previous sections is now $\{(1,2), (1,3), (4,2), (4,3)\}$, corresponding to the 4 "pure" service assignments. Here and for the examples to follow, it is easy to prove that the value functions and related PDE are the same as those one would obtain by allowing appropriate convex combinations of the controls, which can be interpreted as a splitting of effort between different queues.

This problem can be included in the general model of Sections 2 and 4 via an appropriate choice of the $L$ functions. The velocities away from all boundaries for the 4 controls are

$$\bar{\beta}_{(1,2)} = (\lambda - \mu_1, \mu_1 - \mu_2, \mu_2, 0), \qquad \bar{\beta}_{(1,3)} = (\lambda - \mu_1, \mu_1, -\mu_3, \mu_3),$$

$$\bar{\beta}_{(4,2)} = (\lambda, -\mu_2, \mu_2, -\mu_4), \qquad \bar{\beta}_{(4,3)} = (\lambda, 0, -\mu_3, \mu_3 - \mu_4).$$

If $u = (s_A, s_B)$, then we define

$$L_u(\beta) \doteq \begin{cases} -1 & \beta = \bar{\beta}_{(s_A, s_B)} \\ \infty & \text{else.} \end{cases}$$

Note that since the cost effectively constrains $\beta$ once the minimizing control $u$ is chosen, this problem is not actually a game. With this choice the total cost to be minimized is the time to hit the origin, and so $V(x)$ gives the minimal time to reach 0, given the initial state $x$. Even this simple problem is interesting, since the existence of a finite solution to the PDE of Section 4 indicates the existence of a globally stabilizing policy, and in addition suggests a policy that is optimal with regard to the design criteria (returning a perturbed state to zero as quickly as possible).

The PDE in this case takes a very simple form. We have

$$H_{(s_A, s_B)}(\alpha) = 1 + \langle \alpha, \bar{\beta}_{(s_A, s_B)} \rangle,$$

and so $H(\alpha) = \wedge_{(s_A, s_B) \in \mathcal{U}} H_{(s_A, s_B)}(\alpha)$ is a concave function, with polyhedral level sets.

**Example 2: Robust Stabilization.** We next consider an extension of the previous example. In this problem the rates $\lambda, \mu_1, \ldots, \mu_4$ will not be fixed, and so the problem will be a genuine game.

We assume that the rates are restricted to intervals as follows:

$$\lambda \in [a, A] \doteq B_0, \quad \mu_i \in [m_i, M_i] \doteq B_i, i = 1, \ldots, 4.$$

As in the previous example, away from the boundaries there is a well defined velocity associated with each choice of the rates. We wish to consider a game in which the maximizing player will choose the rates, subject to these constraints, in such a way as to delay hitting the origin for as long as possible. Such a game can be obtained by choosing the functions $L_u(\beta)$ to be $-1$ if $\beta$ is in the range of velocities allowed by the service assignment $(s_A, s_B)$, and $\infty$ otherwise. For example

$$L_{(1,2)}(\beta) = \begin{cases} -1 & \text{if } \beta = (\lambda - \mu_1, \mu_1 - \mu_2, \mu_2, 0) \text{ for some } \lambda \in B_0, \mu_1 \in B_1, \mu_2 \in B_2, \\ \infty & \text{otherwise} \end{cases}$$

Note that Condition 4.1 holds for these functions.

The Hamiltonians associated with each assignment are the Legendre transforms of the $L$ functions. We adopt the convention that the operations $\vee$ and $\wedge$ preceed addition and subtraction, and succeed multiplication. Then

$$H_{(1,2)}(\alpha) = c + a\alpha_1 \vee A\alpha_1 + m_1(\alpha_2 - \alpha_1) \vee M_1(\alpha_2 - \alpha_1) + m_2(\alpha_3 - \alpha_2) \vee M_2(\alpha_3 - \alpha_2)$$

$$H_{(1,3)}(\alpha) = c + a\alpha_1 \vee A\alpha_1 + m_1(\alpha_2 - \alpha_1) \vee M_1(\alpha_2 - \alpha_1) + m_3(\alpha_4 - \alpha_3) \vee M_3(\alpha_4 - \alpha_3)$$

$$H_{(4,2)}(\alpha) = c + a\alpha_1 \vee A\alpha_1 + m_2(\alpha_3 - \alpha_2) \vee M_2(\alpha_3 - \alpha_2) + m_4(-\alpha_4) \vee M_4(-\alpha_4)$$

$$H_{(4,3)}(\alpha) = c + a\alpha_1 \vee A\alpha_1 + m_3(\alpha_4 - \alpha_3) \vee M_3(\alpha_4 - \alpha_3) + m_4(-\alpha_4) \vee M_4(-\alpha_4).$$

Hence

$$\begin{aligned} H(\alpha) &= \wedge_{(s_A, s_B) \in \mathcal{U}} H_{(s_A, s_B)}(\alpha) \\ &= c + a\alpha_1 \vee A\alpha_1 + [m_1(\alpha_2 - \alpha_1) \vee M_1(\alpha_2 - \alpha_1)] \wedge [m_4(-\alpha_4) \vee M_4(-\alpha_4)] \\ &\quad + [m_2(\alpha_3 - \alpha_2) \vee M_2(\alpha_3 - \alpha_2)] \wedge [m_3(\alpha_4 - \alpha_3) \vee M_3(\alpha_4 - \alpha_3)]. \end{aligned}$$

The assumptions of Theorem 4.3 apply. Hence a function $V$ satisfying (4.6)–(4.8) will provide the optimal time till the origin is reached. In particular, suppose one can identify a set of optimizing service controls. With these controls, the time till the origin is reached will be bounded by $V(x)$ for any rates $\lambda$ and $\mu_i$ that stay within $B_0$ and $B_i$, $i = 1, \ldots, 4$. A more detailed discussion of the sense in which such a service policy is optimally robust is given in [1].

**Example 3: A Risk-sensitive Problem.** As our final example we consider a problem in risk-sensitive control. Consider a set of nominal constant rates $\bar{\lambda}, \bar{\mu}_1, \ldots, \bar{\mu}_4$ and a parameter $c > 0$. For $z \in \mathbb{R}$ define

$$\ell(z) = \begin{cases} z \log z - z + 1 & \text{if } z \geq 0 \\ \infty & \text{otherwise}, \end{cases}$$

with the convention $0 \log 0 = 0$. As in the previous example, the rates of the system will be chosen by the maximizing player. In contrast to that example, the constraints will not be "hard", although there will be a penalty for any deviation from the nominal rates. The form of the penalty is the one that would be associated with a risk-sensitive control problem, and this point is discussed below. For example, if the service assignment is $(s_A, s_B) = (1, 3)$ then

$$L_{(1,3)}(\beta) = -c + \inf \left\{ \bar{\lambda}\ell(\lambda/\bar{\lambda}) + \bar{\mu}_1\ell(\mu_1/\bar{\mu}_1) + \bar{\mu}_3\ell(\mu_3/\bar{\mu}_3) : \beta = (\lambda - \mu_1, \mu_1, -\mu_3, \mu_3) \right\}.$$

The obvious analogous definitions apply in the other cases. With $h(a)$ defined to be the Legendre transform $(e^a - 1)$ of $\ell$, one can calculate

$$H(\alpha) = c + \bar{\lambda}h(\alpha_1) + \bar{\mu}_1 h(\alpha_2 - \alpha_1) \wedge \bar{\mu}_4 h(-\alpha_4) + \bar{\mu}_2 h(\alpha_3 - \alpha_2) \wedge \bar{\mu}_3 h(\alpha_4 - \alpha_3).$$

The structure of the cost reflects a centering around the nominal system as measured via the function $\ell$. This function is significant from the point of view of large deviations theory, and the game structure allows for a link to related questions in stochastic risk sensitive control theory. Consider a stochastic model in which inter-arrival and service times are exponential random variables with rates $\bar{\lambda}, \bar{\mu}_1, \ldots, \bar{\mu}_4$. Let the time be rescaled by the proportion $n$ and the space by $n^{-1}$, and denote by $\sigma_n$ the (random) time till the origin is reached starting at the fixed point $x$ in the rescaled space. If a service control is applied so as to minimize the cost $\log E e^{nc\sigma_n}$, then the limit of $n^{-1} \log E e^{nc\sigma_n}$ as $n \to \infty$ is formally given by the value $V(x)$ of the game above. As in Example 2, one can prove a robust optimality for a service policy that is designed in this way. This issue is discussed further in [1]